\documentclass[11pt]{amsart}
\usepackage{amsmath,amsthm, amssymb, amsfonts,amscd, epsfig}

\usepackage[mathscr]{eucal}

\theoremstyle{plain}

\newtheorem{thm}{Theorem}[section]
\newtheorem{lem}[thm]{Lemma}
\newtheorem{pro}[thm]{Proposition}
\newtheorem{cor}[thm]{Corollary}
\theoremstyle{definition}
\newtheorem{Def}[thm]{Definition}
\theoremstyle{remark}
\newtheorem{rem}[thm]{Remark}
\newtheorem{exa}[thm]{Example}

\newcommand{\sat}{\mathrm{sat}}
\newcommand{\codim}{\mathrm{codim}}
\newcommand{\depth}{\mathrm{depth}}
\newcommand{\mc}[1]{\mathcal{#1}}

\newcommand{\mbf}[1]{\mathbf{#1}}

\newcommand{\mr}[1]{\mathrm{#1}}

\newcommand{\Ga}{\Gamma}

\newcommand{\Kos}{\mathrm{Kos}}
\newcommand{\Syz}{\mathrm{Syz}}
\newcommand{\Supp}{{\mathrm{Supp} }}

\newcommand{\V}{{\mathbb V}}
\newcommand{\Z}{{\mathbb Z}}

\newcommand{\C}{{\mathbb C}}

\newcommand{\proj}[1]{{\mathbf P}^{#1}}

\begin{document}

\title[Koszul syzygies]{Curvilinear Base Points, 
Local Complete Intersection and Koszul Syzygies in biprojective spaces}

\author{J. William Hoffman and Hao Hao Wang}
\address{Department of Mathematics\\
              Louisiana State University\\
              Baton Rouge, Louisiana 70803}

\email{hoffman@math.lsu.edu, wang\_h@math.lsu.edu}

\begin{abstract}
We prove analogs of results of Cox and Cox/Schenck  
on the structure of certain ideals in the bigraded 
polynomial ring $k[s, u; t, v]$.
\end{abstract}

\thanks{We would like to thank William Adkins
and David Cox for numerous discussions and
suggestions.}

\subjclass{Primary: 14Q10 , Secondary: 13D02, 14Q05}
\keywords{base points, local 
complete intersection, syzygy, saturation, projective space}

\maketitle

\section{Introduction}
\label{s:intro}
Let $X$ be a smooth algebraic surface over a field $k$. 
Let 
$\mc{I} \subset \mc{O}_X$ be a 
coherent sheaf of ideals.  Suppose that $Z=\V(\mc{I})$ is a finite set, 
an assumption that holds throughout this paper.  
Often, we call the points of $Z$ base points for $\mc{I}$.  The reason for the 
terminology is that in our examples $\mc{I}$ will arise from an ideal $I$ in 
a polynomial ring, and we think of generators of $I$ as providing a 
rational map $X \to \proj{N}$. 
Then $Z$ will be the base point locus of this map.  

\begin{Def}
\label{D:curvi}
We say a coherent ideal sheaf $\mc{I} \subset \mc{O}_X$ is 
\textit{curvilinear} at $p\in Z$, if in the completion 
$\widehat{\mc{O}}_{X,p}$, there are analytical coordinates 
$u, v$ so $\widehat{\mc{O}}_{X,p}=k[[u,v]]$ such that 
$\widehat{\mc{I}}_p=\langle u, v^k \rangle$, for some $k \geq 1.$
\end{Def}

\begin{rem}
\label{R:curvi}
In case $X=\proj{2}$, $\mc{I}$ will come from a homogeneous ideal 
$I\subset k[x, y, z]$, and we will say $I$ is curvilinear at $p$ 
if the sheaf it generates is curvilinear. 
Similarly, if $X=\proj{1}\times \proj{1}$, $\mc{I}$ will come from a 
bihomogeneous ideal $I\subset k[s,u,t,v]$.
We will also use the terminology that the base point $p\in Z=\V(\mc{I})$ is 
curvilinear. 
\end{rem}
\begin{Def}
\label{D:LCI}
The base points are \textit{local complete intersection} 
(LCI) if for every point $p\in Z$, $I\mc{O}_{X,p}$ is a complete 
intersection ideal, ie., generated by a regular sequence
in $\mc{O}_{X,p}$. 
\end{Def}

\begin{lem}
If $I$ is curvilinear at $p$, then $I$ is LCI at $p$.
\end{lem}
\begin{proof}
$I\widehat{\mc{O}}_p$ is generated by the regular sequence $u, v^k$, so $I \widehat{\mc{O}}_p$ is LCI.  But $\mc{O}_p/I\mc{O}_p \cong \widehat{\mc{O}}_p/I\widehat{\mc{O}}_p$, because the left hand side is Artinian.  Hence $I\mc{O}_p$ is LCI because $\mc{O}_p$ is a regular ring.
\end{proof}

Recall that given elements $r_1, \cdots, r_n$ in any commutative ring $R$, 
a syzygy $(a_1, \cdots, a_n)$ is a relation 
$a_1r_1+\cdots+a_nr_n=0$ for $a_i \in R$.  
A Koszul syzygy is one of the form $(r_j)r_i+(-r_i)r_j=0$ for 
$i \neq j$.  Let  $\Syz(r_1, \cdots, r_n)$ be the submodule of
 $R^n$ generated by the syzygies.  Let 
$\Kos(r_1, \cdots, r_n)\subset \Syz(r_1, \cdots, r_n)$ be the  
submodule generated by the Koszul syzygies. 
Consider the special case $R=k[x, y, z]$ and $I$ is generated by 
homogeneous forms with $Z=\V(I)\subset \proj{2}$ a finite set of points. 
These syzygy modules are graded and there is no 
loss in generality in considering syzygies that are
homogeneous.
\begin{Def}
\label{D:base}
A syzygy $(a_1, \cdots, a_n)\in \Syz(r_1, \cdots, r_n)$ 
\textit{vanishes at all base points of }$I$ if, for each 
$i$, $a_i\in I^{\sat}$ 
where $I^{\sat}=\{r\in R: \langle x, y, z \rangle ^k r 
\subset I \text{ for some } k \}$.  
If $a_i$ is homogeneous of degree $d_i$, this is equivalent 
to say that $a_i$ belongs to $I\mc{O}_{\proj{2},p}(d_i)$ 
for all $p\in Z$, for all $i$, where $\mc{O}_{\proj{2}}(d)$ is the 
standard invertible sheaf, $d \in \Z$. This is also equivalent to say
that each $a_i$ is a global section of  
$I\mc{O}_{\proj{2}}(d_i) = \mc{I} (d_i )$. 
\end{Def}

The equivalence of the two conditions in the above definition follows 
from the facts
\[ 
I^{\sat}_{d_i}=H^0(\proj{2}, I\mc{O}_{\proj{2}}(d_i)) 
\text{ and } I \mc{O}_{\proj{2}, p}(d_i)=\mc{O}_{\proj{2},p}(d_i) 
\text{ for all  } p\notin Z.
\] 
Now consider an ideal $I=\langle a,b,c \rangle \subset R$ where 
$a,b,c$ are homogeneous polynomials of degree $m$. 
Cox has proved following theorem:

\begin{thm}
\label{T:C}(\cite[Theorem 5.5]{C})
Let $k$ be algebraically closed. If all the base points 
of $I$ are curvilinear, then every 
syzygy on $a,b,c$ which vanishes on $Z$ is Koszul syzygy.  
\end{thm}

Cox proves this under the assumption that 
$k= \C$, but the arguments are more general. 
We will prove a version of this theorem for the bigraded case 
in Section \ref{s:cl}.
Cox and Schenck \cite{CS} have proved the following result which 
strengthens Theorem \ref{T:C}.

\begin{thm} (\cite[Theorem 1.7]{CS})
\label{T:CS} 
If $I=\langle f_1, f_2, f_3 \rangle \subset R=k[x, y, z]$ where 
$f_i$ is a homogeneous polynomial, and if $I$ has codimension $2$,  
then the module of syzygies of $f_1, f_2, f_3$ vanishing at 
$Z = \mathbb{V}(I)$ is generated by the Koszul 
syzygies if and only if $I$ is a local complete intersection. 
\end{thm}

We will extend this result to the bigraded case in Section \ref{s:lci}.
The proofs of theorems (\ref{T:C}) and (\ref{T:CS}) are quite 
different; the former is geometric, the latter algebraic. 
We present arguments that parallel the proofs of the above 
theorems. Even though the first follows logically from the second, 
we present both arguments, because they are so different. 
Also, there are some significant differences from the case considered
by Cox and Schenck. In fact we can prove only that those syzygies
vanishing at the base points, {\it and of sufficiently high degree}, 
are Koszul syzygies. We present an example to show that, in the 
bigraded case, not every syzygy vanishing at the base points
is Koszul. The key difference between their situation and 
ours is that, for any invertible sheaf $\mc{L}$ on $\proj{2}$, 
we have $H^1 (\proj{2},\, \mc{L}) = 0$, whereas this vanishing
does not hold on $\proj{1}\times \proj{1}$. The results in this 
paper are used in an essential way in Wang's thesis, \cite{Wang}.

\section{Curvilinear base points and Koszul syzygies}
\label{s:cl}
In this section, $k$ is an algebraically closed field.
Let $R=k[s,u,t,v]$ be the bigraded coordinate ring of
$\proj{1}\times \proj{1}$, where $s,u$ have bidegree $(1,0)$ 
and $t,v$ have bidegree $(0,1)$.
Let $\mbf{m}=\langle st, sv, ut,uv \rangle$, the irrelevant ideal in $R$.
 Consider an ideal 
$I=\langle a,b,c \rangle \subset R$, where $a,b,c$ are homogeneous of
the same degree $(m,n)$. Let $Z=\V(a,b,c) \subset \proj{1}\times \proj{1}$,
which is assumed to be a finite set. 
Let $r_1, \cdots, r_n \in R$ be bihomogeneous elements, 
and consider syzygies as before.  
All our syzygy modules will be bigraded, and syzygies themselves
will be taken to be bihomogeneous. 
\begin{Def}
\label{D:base2}
A syzygy $(a_1, \cdots, a_n)\in \Syz(r_1, \cdots, r_n)$ 
\textit{vanishes at all base points of} $I$ if, for each 
$i$, $a_i \in I^{\mbf{m}\sat}=\{r\in R: \mbf{m} ^k r \subset I, 
\text{ for some } k\}$. If $a_i$ homogeneous of bidegree 
$(d_i, d'_i)$, this is equivalent to say that 
$a_i\in I\mc{O}_{\proj{1}\times \proj{1},p}(d_i, d'_i)$, for all $p\in Z$, 
where $\mc{O}_{\proj{1}\times \proj{1}}(d, d')$ is the 
standard invertible sheaf, $d, d' \in \Z$. This is also equivalent to say
that each $a_i$ is a global section of  
$I\mc{O}_{\proj{1}\times \proj{1}}(d_i, d_i ') = \mc{I} (d_i , d_i ')$. 
\end{Def} 
The equivalence follows as before, noting that 
$I^{\mbf{m}\sat}_{d_i,d'_i}=H^0(\proj{1}\times \proj{1}, 
I\mc{O}_{\proj{1}\times \proj{1}}(d_i,d_i'))$ and 
$I\mc{O}_{\proj{1}\times \proj{1},p}(d_i,d'_i)=
\mc{O}_{\proj{1}\times \proj{1},p}(d_i,d'_i)$ for all $p\notin Z$.
We shall be concerned with the syzygies $(A,B,C)$ on the generators of $I=(a,b,c)$.
\begin{rem}\label{R:puredegree}
Consider the exact sequence:
\[
\begin{CD} 0@>>> \Syz(a,b,c) @>>> \bigoplus_{i=1}^3 R(-d_i, -d'_i) 
@>(a,b,c)>>I@>>>0 \end{CD}.\]
$\Syz(a,b,c)$ is a bigraded submodule of 
$\bigoplus_{i=1}^3 R(-d_i, -d'_i)$, and we say 
$(A,B,C)\in \Syz(a,b,c)_{k,l}$ has bidegree $(k,l)$.  
Note that, as polynomials, the bidegree of $A$ is $(k-d_1, l-d'_1)$, 
since $A \in R(-d_1, -d'_1)_{k,l}$, etc. 
In other words, the polynomial expression $Aa+Bb+Cc$ has bidegree $(k,l).$
If all $(d_i, d'_i)$ are equal to a fixed pair $(m,n)$, then a 
bihomogeneous syzygy $(A,B,C)$ of bidegree $(k,l)$ will have $A,B,C$ all 
bihomogeneous of same bidegree $(k-m, l-n)$.  Some authors call this a 
syzygy of bidegree $(k-m,l-n)$. 
We will call $(A,B,C)$ a syzygy of \textit{pure degree} $(k-m,l-n)$. 
\end{rem}

\begin{pro}\label{P:vankoz}
Let $a,b,c \in R$ be bihomogeneous polynomials of bidegree $(m,n)$ and 
suppose $\V(a,b,c)\subset \proj{1}\times \proj{1}$ is finite. Suppose that all the base points of
$a,b,c$ are curvilinear. Then every syzygy $Aa+Bb+Cc=0$ with
pure bidegree $(k,l)$, where $(k-2m+1)(l-2n+1) \geq 0$, and 
vanishing at the base points is a Koszul syzygy. That is, there exists 
$h_1,h_2,h_3$ of bidegree $(k-m,l-n)$ such that
\begin{align*}
A&=h_1 c+h_2 b\\
B&=-h_2 a+h_3 c\\
C&=-h_1 a-h_3 b
\end{align*}
\end{pro}

Note: when $\V(a,b,c)=\emptyset$, this is proved in \cite{CGZ}. The
analogous result for $\phi : \proj{2} \longrightarrow \proj{3}$, 
namely Theorem \ref{T:C}, is proved in \cite{C}, but with no 
restriction at all on the degree of the syzygies. Note also that 
the Koszul syzygies themselves have pure degree $(m, n)$ and
satisfy the conditions of the proposition. 
 
\begin{proof}
In \cite{C}, using the theory of toric varieties, Cox constructed a
projective birational map (blowing up) 
$\pi: X \to Q=\proj{1}\times \proj{1}$ with the following properties:
\begin{itemize}
\item[1.] $X$ is a smooth projective surface;
\item[2.] $I'_Z=I\mc{O}_X=\mc{O}_X(-E)$ for a divisor $E$ on $X$;
\item[3.] The canonical class 
$K_X=\pi^{*} K_Q+E$ where $K_Q=\mc{O}_Q(-2,-2)$ is the canonical class
of $Q$. 
\end{itemize}
 We have $a,b,c \in \Ga(Q, \mc{I}_{Z}\otimes
\mc{O}_{Q}(m,n)) \subset R_{m,n},$ and we get their proper
transform 
$$\tilde{a},\tilde{b},\tilde{c} \in \Ga(X, \mc{I}_{Z}^{'}
\otimes \pi^{*} \mc{O}_{Q}(m,n)).$$
Let $L_1,L_2$ be generators of
$\mr{Pic}(Q)$, i.e. generating  lines on $Q$ which are chosen to
miss $Z$. By abuse of notation, we let $L_1,L_2 \subset X$ be the
inverse images on $X$. Hence 
$$\tilde{a},\tilde{b},\tilde{c} \in
\Ga(X,\mc{O}_{X}(mL_1+nL_2-E)).$$
These are without base points, i.e., 
\[ 
\begin{CD} 
\mc{F} :=  \mc{O}_{X}(E-mL_1-nL_2)^3
@>{(\tilde{a},\tilde{b},\tilde{c})}>> \mc{O}_{X} \longrightarrow 0
\end{CD} 
\]
is onto. Build the Koszul complex on this, i.e.
\begin{equation}
\label{E:koz}
0 \rightarrow \wedge^{3}\mc{F} \rightarrow \wedge ^{2} \mc{F}
\rightarrow \mc{F}\rightarrow \mc{O}_X \rightarrow 0 
\end{equation}
with $\wedge ^3 \mc{F}= \mc{O}_{X}(3E-3mL_1-3nL_2)$, and $\wedge ^2
\mc{F} = \mc{O}_{X}(2E-2mL_1-2nL_2)^3$. This is exact on $X$,
because $\tilde{a}, \tilde{b}, \tilde{c}$ have no base points
\cite[Prop. 17.14(b), p.439]{E}.
\par
Let $A,B,C$ be a syzygy of pure bidegree $(k,l)$ which vanishes at base
points, i.e., 
$$A,B,C \in \Ga(Q, \mc{I}_{Z}\otimes \mc{O}_Q(k,l)).$$
Their proper transforms $\tilde{A},\tilde{B},\tilde{C} \in
\Ga(X,\mc{O}_{X}(kL_1+lL_2-E))$.  Tensor the Koszul sequence 
(\ref{E:koz}) with
$\mc{O}_{X}((k+m)L_1+(l+m)L_2-2E)$ and break the resulting
sequence of 4 terms into 2 exact sequences
\[
0 \rightarrow \mc{O}_{X}((k-2m)L_1+(l-2n)L_2+E) \rightarrow
\mc{O}_{X}((k-m)L_1+(l-n)L_2)^3 \rightarrow \mc{K} \rightarrow 0
\]
\[
0 \rightarrow \mc{K} \rightarrow \mc{O}_{X}(kL_1+lL_2-E)^3
\rightarrow \mc{O}_{X}((k+m)L_1+(l+n)L_2-2E) \rightarrow 0 
\]
From the exact cohomology sequence, we see that if
\[ 
H^1(X,\mc{O}_{X}((k-2m)L_1+(l-2n)L_2+E))=0
\]
then the sequence
\[
H^0(X,\mc{O}_{X}((k-m)L_1+(l-n)L_2))^3 \rightarrow
H^0(X,\mc{O}_{X}(kL_1+lL_2-E))^3 
\]
\[
\begin{CD} 
@>{(\tilde{a},\tilde{b},\tilde{c})}>>
H^0(X,\mc{O}_{X}((k+m)L_1+(l+n)L_2-2E)) 
\end{CD}
\] 
is exact in the middle.  But because 
$\pi_{*} \mc{O}_{X}=\mc{O}_{Q}$ and
$ R^i \pi_* \mc{O}_X= 0$ if $i > 0$ (see below Lemma), we have that
\[H^0(X,\mc{O}_{X}((k-m)L_1+(l-n)L_2)) \cong
H^0(Q,\mc{O}_{Q}(k-m,l-n)). \] This follows from the 
Leray spectral sequence. Note that 
$\mc{O}_{X}(kL_1+lL_2)=\pi^{*} \mc{O}_{Q}(k,l)$. Then $\pi_{*}
\mc{O}_{X}(kL_1+lL_2)= \mc{O}_{Q}(k,l) \otimes \pi_{*} \mc{O}_{X}$
by the projection formula \cite[Ch. III, \S8, Ex. 8.3, p 253]{H}.

We can consider the following diagram:
\[
\begin{CD}
 H^0(Q, \mc{O}_Q(k-m,l-n))^3 @>\lambda>> 
         H^0(X, \mc{O}_X((k-m)L_1+(l-n)L_2))^3\\
@V\gamma VV  @V\alpha VV\\
H^0(Q, I_Z\otimes \mc{O}_Q(k,l))^3 @>\mu >> 
H^0(X, \mc{O}_X(kL_1+lL_2-E))^3\\
@V\delta VV  @V\beta VV\\ 
H^0(Q, I^2_Z\otimes \mc{O}_Q(k+m,l+n)) @>\nu >>
H^0(X,\mc{O}_X((k+m)L_1+(l+n)L_2-2E))
\end{CD}
\]
where $\delta = (a, b, c)$ and  
$\beta=(\tilde{a}, \tilde{b}, \tilde{c})$.
We have seen that $\lambda$ is an isomorphism, 
and clearly, $\mu$ and $\nu$ are injective.
Now an element $(A, B, C)$ from the middle
term of the first column in the kernel of $\delta$ 
is precisely a syzygy vanishing at the base points. 
$Aa+Bb+Cc=0$ implies that 
$(\tilde{A}, \tilde{B}, \tilde{C})$ is in the 
image of $\alpha$, since the right sequence is exact in the middle, 
as we have shown. 
Thus, $A,B,C$ is in the image of $\gamma$, since the diagram commutes. 
But this says exactly that $(A, B, C)$ is a Koszul syzygy, 
as was to be shown.
\par
It remains to check that
$H^1(X,\mc{O}_{X}((k-2m)L_1+(l-2n)L_2+E))=0$. By Serre duality,
this is dual to
\[H^1(X,\mc{O}_{X}(-(k-2m)L_1-(l-2n)L_2-E)\otimes \mc{K}_{X})\]
 but $\mc{K}_{X} \cong \pi^{*} \mc{K}_{Q}+E$, and 
$\mc{K}_{Q} \cong \mc{O}_{Q}(-2,-2)$.  Hence, the
above is
\[ 
H^1(X,\mc{O}_{X}((2m-k-2)L_1+(2n-l-2)L_2))=
H^1(X,\pi^* \mc{O}_{Q}(2m-k-2,2n-l-2)). 
\]
From the Leray spectral sequence, the Lemma below, and
the projection formula, we get
\[ 
H^1(X, \pi^* \mc{O}_{Q}(2m-k-2,2n-l-2)) 
\cong H^1(Q,\mc{O}_{Q}(2m-k-2,2n-l-2)).
\]
Now apply the K\"unneth formula, \cite{SW}:
\begin{align*}
& H^1 (\proj{1}\times \proj{1}, \mc{O}_{\proj{1} \times \proj{1}}(2m-k-2,2n-l-2)) \\
&=\bigoplus_{i+j=1} H^{i} (\proj{1}, \mc{O}(2m-k-2)) \otimes 
H^j (\proj{1},\mc{O}(2n-l-2))
\end{align*}  
By Serre's computation of the cohomology of 
projective space and this last one clearly is 0, when
$k=2m-1$, $\forall l$, 
or $l=2n-1$, $\forall k$, or $k >2m-1$, $l>2n-1$, or $k<2m-1$, $l<2n-1$.
i.e. $(k-2m+1)(l-2n+1)\geq 0$.
\end{proof}

\begin{lem}
\label{L:vanish}
\[
\pi_* \mc{O}_X= \mc{O}_Q,\ \  R^i \pi_* \mc{O}_X= 0 \text{ if} \ \ i > 0 
\]
\end{lem}
\begin{proof}
Since $\pi$ is a proper birational morphism of smooth
surfaces, Zariski's Factorization Theorem 
\cite[Corollary 5.4, Page 411]{H} states that we may factor
$\pi$ into a succession of maps
$\pi _i : X_i \to X_{i-1}$ where $X_0 = Q$, $X_n = X$, 
and each $\pi _i$ is a blowing up of one point (monoidal
transformation). The result now follows by induction on $n$, 
using the spectral sequence of a composite functor, and the 
fact, shown by direct computation, 
that it is true for a monoidal transformation. In fact, this
result is a well-known fact for any succession of blowing-ups
with smooth centers.   
\end{proof}
\begin{rem}
\label{R:general}
The reader may wonder where the hypothesis of curvilinear
base point has come into the above argument. It occurs in the 
existence of the blowing-up $\pi$ with the properties enumerated
at the beginning of the proof. Starting with {\it any} ideal
$I$ one may construct a blowing-up to a smooth $X$ 
with the property that
$I\mc{O}_X=\mc{O}_X(-E)$ for a divisor $E$ on $X$.  
This is a general fact for varieties in any dimension, 
at least in characteristic $0$, by Hironaka's theorems, 
and one may even assume $E$ to have normal crossings. In 
our case, for surfaces, Zariski's theorem utilized above
makes it clear that 
$K_X=\pi^{*} K_Q+E'$ for some divisor $E'$ supported
in the same set as $E$, but the key fact in the curvilinear
case is that $E = E'$, an exact equality including multiplicities.
In view of the fact that the conclusion of the theorem holds
for $I$ being local complete intersection, it might be of 
interest to extend the above proof to that case. This requires
some analysis of the divisor $E'$ and it's relation to 
$E$ in that situation.  
\end{rem}

\begin{exa}
The following example shows that not all syzygies vanishing at the base
points are Koszul syzygies: Let 
$I=\langle u^2tv, u^2t^2+suv^2, s^2tv \rangle$.  
The base points of ideal
$I$ are 
$$p=(0:1;0:1), \ p'=(1:0;1:0),$$  and 
$$I_p=\langle s,t \rangle, \ I_{p'}=\langle v, u^2 \rangle.$$  
These are curvilinear. The (pure) degrees of the 
generators are $(2, 2)$ so our theorem claims
that all syzygies vanishing at the base points and either of degrees
$\ge (3, 3)$ or $\le (3, 3)$ are generated by Koszul syzygies.   
Using Singular \cite{GPS01}, we found generators for the syzygy module vanishing 
at the base points as the rows of the following matrix:
\[ \left[\begin{matrix}
 u^2t^2+suv^2 & -u^2tv   & 0 \\
s^2tv         & 0        &-u^2tv \\
0             & s^2tv    &-u^2t^2-suv^2 \\
ut^3v+stv^3   &-ut^2v^2  & 0 \\
st^3v         &-st^2v^2  &utv^3\\
s^2t^2v       &-s^2tv^2   &suv^3\\
s^3uv         &0         &-su^3v\\
s^2u^2t       &0         &-u^4t\\
s^3u^2        &0         &-su^4
\end{matrix}\right].\]
The Koszul syzygies are generated by the first three rows.
The sixth row is of degree $(2,3)$ and it is generated by the Koszul syzygies since 
\[(s^2t^2v, -s^2tv^2,suv^3)=t(s^2tv,0,-u^2tv)-v(0,s^2tv,-u^2t^2-suv^2).\]
The syzygies
\[(ut^3v+stv^3,-ut^2v^2,0), (st^3v, -st^2v^2,utv^3), (s^3uv, 0, -su^3v), (s^2u^2t, 0,-u^4t), (s^3u^2, 0, -su^4)\] are of degree $(1,4), (1, 4), (4,1), (4, 1), (5,0)$.  
They are not Koszul, 
because their degrees forbid this.
But we can look at multiples of say $ \mbf{v} = (st^3v, -st^2v^2, utv^3)$, 
which is of degree $(1,4)$.  
Since
\[s(st^3v, -st^2v^2, utv^3)=t^2(s^2tv, 0,-u^2tv)-tv(0, s^2tv,-u^2t^2-suv^2).\]
we see that $s^2 \mbf{v}$ and $su \mbf{v}$, of degrees $(3, 4)$, are in the module of Koszul syzygies, 
as is claimed by our theorem. 
 The syzygy $u\mbf{v}$, of degree $(2,4)$, is not generated by Koszul syzygies, 
but 
\[u^2(st^3v, -st^2v^2, utv^3)=stv(u^2t^2+suv^2, -u^2tv,0)-uv^2(s^2tv,0,-u^2tv),\]
so $u^2\mbf{v}$, of bidegree $(3,4)$, is in the module of Koszul syzygies. 
\end{exa}

\section{Local complete intersection and Koszul syzygies} 
\label{s:lci}
Now assume that $k$ is any field and consider an ideal 
$I=\langle f_1, f_2, f_3\rangle \subset 
R = k[s, u, t, v]$, where $f_i$ is 
bihomogeneous of degree $(d_i, d_i')$.  The $f_i$ form 
a regular sequence in $R$ if and only if the following 
Koszul complex is exact.
\begin{equation}
\label{E:exact}
 \begin{CD} 0 \rightarrow R(-\sum_{i=1}^{3}d_i,
-\sum_{i=1}^{3}d_i')
 @>{\left[ \begin{matrix} f_3 \\ -f_2 \\ f_1 \end{matrix} \right]} >>
 \bigoplus_{i<j} R(-d_i-d_j, -d_i'-d_j')
\end{CD} 
\end{equation}
 \[ 
\begin{CD} @>{\left[ 
\begin{matrix} f_2 & f_3 & 0 \\
 -f_1 & 0 & f_3 \\ 0 & -f_1 & -f_2 \end{matrix} 
\right]} >>
 \bigoplus_{i=1}^{3} R(-d_i, -d_i') @>{ [f_1 \ \  f_2 \ \  f_3
 ] } >> I \rightarrow 0
\end{CD}.
\]
We will discuss the situation when $Z=\V(I) \subset \proj{1}
\times \proj{1}$ is a zero-dimensional subscheme, thus $I=\langle
f_1, f_2, f_3 \rangle$ has codimension two in $R$.  We call $Z$
the base point locus of $f_1, f_2, f_3$. If $I$ has codimension
two, $f_1, f_2, f_3$ will no longer be a regular sequence: since
$R$ is Cohen-Macaulay, the depth of $I$ equals the codimension 
of $I$, so that a maximal $R$-sequence in $I$ has length two.

\begin{lem}
If $I=\langle f_1, f_2, f_3 \rangle$ has codimension two in $R$,
then the sequence \eqref{E:exact} is exact except at $\bigoplus_{i=1}^{3} R(-d_i,
-d_i')$.  In particular, the Koszul complex of $f_1, f_2, f_3$
gives the exact sequences:
\[  0 \rightarrow R(-\sum_{i=1}^{3}d_i,
-\sum_{i=1}^{3}d_i') \rightarrow
 \bigoplus_{i<j} R(-d_i-d_j, -d_i'-d_j')\rightarrow \bigoplus_{i=1}^{3} R(-d_i,
 -d_i')\]
 and
 \[ \bigoplus_{i=1}^{3} R(-d_i, -d_i') \rightarrow I \rightarrow 0
 \]
\end{lem}

\begin{proof}
The exactness of the first sequence follows from
$\depth(I)=\codim(I)=2$ and the Buchsbaum-Eisenbud exactness
theorem (\cite[p. 500]{E}).
\end{proof}

\begin{Def}
A \textit{Koszul syzygy} on $f_1, f_2, f_3$ is an element of the
submodule
\[ K \subset \bigoplus_{i=1}^{3} R(-d_i, -d_i') \] generated by
the columns of the matrix
\[ \left[ \begin{matrix} f_2 & f_3 & 0 \\
 -f_1 & 0 & f_3 \\ 0 & -f_1 & -f_2 \end{matrix} \right] \]
\end{Def}

\begin{cor} If $I=\langle f_1, f_2, f_3 \rangle$ is a codimension
two ideal, then we have an exact sequence 
\begin{equation}\label{E:K}
0 \rightarrow R(-\sum_{i=1}^{3}d_i,
-\sum_{i=1}^{3}d_i') \rightarrow
 \bigoplus_{i<j} R(-d_i-d_j, -d_i'-d_j')\rightarrow K \rightarrow 0 
\end{equation}
\end{cor}
Note: $K$ is a submodule of the syzygy module $S$ defined
by the exact sequence
\[ 0 \rightarrow S \rightarrow \bigoplus_{i=1}^{3} R(-d_i, -d_i')
\rightarrow I \rightarrow 0. \] In fact, $K$ is a proper submodule 
of $S$, because if $K=S$, we would have exactness of the sequence 
\eqref{E:exact}, which is equivalent to the $f_i$ forming a regular
sequence, which we have remarked, is impossible.

Analogous to \cite[page 32]{SR}, we have the following
definition:
\begin{Def}
A bigraded submodule $M$ of a finitely generated free $R-$module
$F$ is \textit{saturated} if
\[ 
M=\{ x \in F | \mbf{m} x \subset M  \}
\]
We define the \textit{saturation} of $M$ to be 
\[
M^{\mbf{m}\sat}=\{x\in F | \mbf{m} ^k x \subset M, \text{ for some } k \} .
\]
where $\mbf{m} = \langle st, sv, ut, uv \rangle$. We will use the
notation $M^{\sat}$ for this, the ideal $\mbf{m}$ being understood. 
\end{Def}

\begin{rem}
$M$ is saturated if and only if $M=M^{\sat}$.
\end{rem}
\begin{proof}
Suppose $M$ is saturated, we will prove that $M=M^{\sat}$.
Since $M \subset M^{\sat}$, we only need to show 
$M^{\sat} \subset M$ when $M$ is saturated.  Let $x \in M^{\sat}$. 
There exists a $k$ such that 
$\mbf{m}^k x =\mbf{m} \mbf{m}^{k-1} x \subset M$. 
Since $M$ is saturated, we have $\mbf{m}^{k-1}x \subset M$.  
We can repeat the process until $\mbf{m} x \subset M$.  
Therefore $x \in M$. Thus, $M^{\sat} \subset M$, and we proved that $M=M^{\sat}$. 
\par
On the other hand, suppose $M=M^{\sat}$,
if $x\in F$ such that $\mbf{m} x \in M$, then $x \in M^{\sat}=M$. 
Therefore, $M$ is saturated.
\end{proof}

\begin{pro}
\label{P:sat}
Let $M$ be a bigraded submodule of a free
$R=k[s, u, t, v]$-module of finite rank $F$. 
Let $\mc{M}$ be the corresponding coherent sheaf on 
$X = \proj{1}\times \proj{1}$. Then 
\[
M^{sat}_{k, l} = H^0 (X, \, \mc{M} (k, l)).
\]
\end{pro}
\begin{proof}
For any finitely generated bigraded $R$-module we have an 
exact sequence (see \cite{HY})
\[
\label{E:local1}
\begin{CD}
0 @>>> H^0_{\mbf{m}}(M) @>>> M @>>>
\displaystyle{\bigoplus _{(a,b) \in \Z ^2}H^0(X,\,  \mc{M}(a,b))} @>>>
H^1_{\mbf{m}}(M) @>>>0.
\end{CD}
\]
We will show that  $H^i_{\mbf{m}}(M)=0$ when $i=0, 1$ if $M$ is saturated.
This is sufficient since both $M$ and $M^{sat}$ generate the same sheaf.
Since  $H^i_{\mbf{m}}(R)=0$, for $i = 0, 1$,
and $M$ is a submodule of a free $R$-module, 
it is clear that we have vanishing for $i=0$, and this does not depend
on $M$ being saturated. The long exact cohomology sequence for
\[
\begin{CD}
0 @>>> M @>>> F @>>> F/M @>>> 0
\end{CD}
\]
has a piece
\[
\begin{CD}
 H^0_{\mbf{m}}(F )@>>> H^0_{\mbf{m}}(F/M) @>>>   
H^1_{\mbf{m}}(M ) @>>>
H^1_{\mbf{m}}(F ) 
\end{CD}.
\]
Since the extreme terms are zero, we get an isomorphism
\[
H^1_{\mbf{m}}(M ) \cong  H^0_{\mbf{m}}(F/M)
\]
But the right-hand side is clearly $M^{sat}/M$, proving
our claim. 
\end{proof}

We will use the following well-known results:

\begin{pro}
\label{P:supp}
Let $R$ be a Noetherian ring, $J \subset R$ an ideal, and
$M$ a finitely
generated $R$-module. The following are equivalent:
\begin{itemize}
\item[1.] ${\rm Ass} (M) \subset \mathbb{V}(J)$.
\item[2.] ${\rm Supp} (M) \subset \mathbb{V}(J)$.
\item[3.] There exists $n \ge 0$ such that $J^n M = 0$. 
\end{itemize}
When this is so, $H^0_J(M)=M$.
\end{pro}
\begin{proof}
The equivalence of (1) and (2) follows from the fact that
both  ${\rm Ass} (M)$ and  ${\rm Supp} (M)$ have the same
minimal elements (\cite[Theorem 1, p. 7]{Serre}). 
Assume (3) holds, and let $\mbf{p}$ 
be  prime ideal in the support of $M$. 
Let $m/s \in M_{\mbf{p}}$.  If 
$\mbf{p}\nsupseteq J$ there would exist 
$x\in J \setminus  \mbf{p}$, and clearly 
$x^n\in J^n \setminus \mbf{p}$ for any $n \ge 1$. 
But then $x^n m = 0$, and this shows that
$m/s = 0$, showing that $\mbf{p}$ cannot be in the support 
of $M$. Conversely, assume (2) then (\cite[Prop 3, p. 5]{Serre}) 
\[
{\rm Supp} (M) = \mathbb{V}({\rm ann} (M)) \subset \mathbb{V}(J)
\]
shows that $\sqrt {J} \subset \sqrt{{\rm ann} (M)}$, from
which (3) follows easily. 

To see the last statement, note that $H^0_J(M)=\{ m\in M:J^k m=0 \text{ for some } k\}$. 
\end{proof}

\begin{rem}
\label{R:support}
If $M$ is a bigraded module over the ring
$R = k[s, u, t, v]$, the elements of   ${\rm Ass} (M)$
are bihomogeneous. Moreover, taking  $J = \mbf{m}$ the
irrelevant ideal, the conditions of the previous 
proposition are easily seen to be equivalent to
\begin{itemize}
\item[4.]
There exists ${m, n}$ such that  $M _{k, l} = 0$ whenever
$k \ge m$ and $l \ge n$, which we abbreviate by writing
$(k, l) \gg (0, 0)$.
\end{itemize} 
\end{rem}

\begin{Def}
The \textit{bigraded Hilbert polynomial} $H(M)$ of a finitely
generated bigraded $R$-module $M$ is the unique polynomial in two 
variables such
that
\[ H(M)(n, n') = \dim_k M_{n, n'} \]
for all $n, n' \gg 0,$ where $M_{n,n'}$ is the bigraded piece of
$M$ in degree $(n,n')$.
\end{Def} 

Note: if $\mc{M}$ is the sheaf of modules associated to $M$, for $n,n' \gg 0$, we have 
\[ H(M)(n, n') = \dim_k M_{n, n'}=\dim H^0(\proj{1} \times \proj{1},\mc{M}(n,n')). \]

\begin{rem}
\label{R:sat}
\[H(M)=H(M^{\sat}).\]
\end{rem}
\begin{proof}
Let $Q=M^{\sat}/M$.  Since $\mbf{m}^r Q =0$ for some $r$, 
we must have $Q_{k,k'}=0$ for $k,k' \gg 0$, as in the previous remark.  
Thus $H(M)=H(M^{\sat})$.
\end{proof}

\begin{lem}
Let $M\subset N \subset  F$ be bigraded submodules where $F$ is free of 
finite type.  If $M$ is saturated,  then $M=N$ if and only if $H(M)=H(N)$.
\end{lem}

\begin{proof}
We will only show that $H(M)=H(N)$ implies that $M=N$. $H(M)=H(N)$ 
says that there exist $n,n'$ such that 
$M_{k, k'}=N_{k,k'}$ for all $k \geq n, k'\geq n'$.  Let $a_{p, p'} \in
N_{p,p'}$ where $p<n$, or $p'<n'$.  We can find an $\alpha$ such 
that $\mbf{m}^{\alpha} a_{p,p'} \subset N_{k,k'}=M_{k,k'}$ 
for some $k\geq n, k'\geq n'$.  
Since $M$ is saturated, $M=M^{\sat}$, thus $a_{p,p'}\in M_{p,p'}$.
\end{proof} 

As before, we denote $S$ as the syzygy module, $K$ as the Koszul syzygies 
and $V$ as the module of syzygies for $f_1, f_2, f_3$ vanishing at the base points. 

\begin{lem}
\label{L:sat}
As submodules of  $\bigoplus_{i=1}^{3} R(-d_i, -d_i')$, 
$V$ is a saturated submodule, and  
$K_{k,k'}=K^{\sat}_{k,k'}$ when $(k-\sum_{i=1}^3d_i+1)(k'-\sum_{i=1}^3 d'_i+1)\geq 0$. 
\end{lem}

\begin{proof}
We know that $K \subset V \subset \bigoplus_{i=1}^{3} R(-d_i,
-d_i')$. 
We first consider $V$. By definition, we
have
\[ 
V=S \cap \bigoplus_{i=1}^{3} I^{\sat}(-d_i, -d_i').
\]
Note that $S=S^{\sat}$, or equivalently, that $S$ is saturated.  
To see this, let $(a,b,c)\in R^3$  such that $\mbf{m}(a,b,c) \subset S$.  
This means that for all $h \in \mbf{m}$ and $h(a,b,c) \in
\Syz(f_1, f_2, f_3)$.  This says that $h(af_1+bf_2+cf_3)=0$. But $R$
has no zero divisors, thus
$(a,b,c)\in S$, which shows that $S$ is saturated. 
Since the intersection of saturated
submodules is saturated, $V$ is saturated.
\par
We will show $K^{\sat}_{k,k'}=K_{k,k'}$ for all  
$(k-\sum_{i=1}^3 d_i+1)(k'-\sum_{i=1}^3 d'_i+1)\geq 0$.  
Let $r=r_{k,k'}\in \bigoplus_{i=1}^{3}
R(-d_i, -d_i')$ with bidegree $(k-\sum_{i=1}^3 d_i+1)(k'-\sum_{i=1}^3 d'_i+1)\geq 0$ and  satisfies $\mbf{m}r
\subset K$.  We will show that $r \in K$.  Let $L=K+Rr$. 
Consider the short exact sequence
\[ 0 \rightarrow K \rightarrow L \rightarrow L/K \rightarrow 0.\]
We get a long exact
sequence in local cohomology
\[ 0 \rightarrow H_{\textbf{m}}^0(K) \rightarrow  H_{\textbf{m}}^0(L)
\rightarrow H_{\textbf{m}}^0(L/K) \rightarrow H_{\textbf{m}}^1(K)
\rightarrow \] 
Since $ K \hookrightarrow \bigoplus_{i=1}^{3} R(-d_i, -d_i')$, 
$H^0_{\mbf{m}}(K)=0$.  Also consider the exact local cohomology 
sequence of the exact sequence \eqref{E:K}
\[ 
\begin{CD}
\bigoplus_{i <j} H^i_{\mbf{m}}(R(-d_i-d_j, -d'_i-d'_j))_{k,k'} 
@>>> H^i_{\mbf{m}}(K)_{k,k'} @>>> H^{i+1}_{\mbf{m}}
(R(-\sum_{i=1}^3d_i, -\sum_{i=1}^3d'_i)) 
\end{CD}.
\]
$H^1_{\mbf{m}}(K)_{k,k'}=0$ if 

\begin{eqnarray}
\label{E:R1} \bigoplus_{i <j} H^1_{\mbf{m}}(R(-d_i-d_j, -d'_i-d'_j))_{k,k'}=0  
\text{ and }\\
\label{E:R2}  H^2_{\mbf{m}}(R(-\sum_{i=1}^3d_i, -\sum_{i=1}^3d'_i))_{k,k'}=0.
\end{eqnarray}

Since we know that $R$ is strongly $(0,0)$-regular (see \cite{HW}), 
Equation \eqref{E:R1} holds for all $k,k'$, and Equation \eqref{E:R2} 
can be written as the following
\[ H^2_{\mbf{m}}(R)_{k-\sum_{i=1}^3d_i, k'-\sum_{i=1}^3d'_i}=H^1(\proj{1}\times \proj{1}, \mc{O}_{\proj{1}\times \proj{1}}(k-\sum_{i=1}^3d_i, k'-\sum_{i=1}^3d'_i))=0.\]
Now apply the K\"unneth formula, \cite{SW}:
\begin{align*}
& H^1 (\proj{1}\times \proj{1}, \mc{O}_{\proj{1} \times \proj{1}}(k-\sum_{i=1}^3d_i, k'-\sum_{i=1}^3d'_i)) \\
&=\bigoplus_{i+j=1} H^{i} (\proj{1}, \mc{O}(k-\sum_{i=1}^3d_i)) \otimes 
H^j (\proj{1},\mc{O}(k'-\sum_{i=1}^3d'_i))
\end{align*}
By Serre's computation of the cohomology of 
projective space and this last one clearly is 0, when
$k=\sum_{i=1}^3d_i-1$, $\forall l$, 
or $l=\sum_{i=1}^3d'_i-1$, $\forall k$, or $k >\sum_{i=1}^3d_i-1$, $l>\sum_{i=1}^3d'_i-1$, or $k >\sum_{i=1}^3d_i-1$, $l>\sum_{i=1}^3d'_i-1$.
i.e. $(k-\sum_{i=1}^3 d_i+1)(k'-\sum_{i=1}^3 d'_i+1)\geq 0$.

Therefore, we have
\[
H^1_{\mbf{m}}(K)_{k,k'}=0, \ \ \text{when } (k-\sum_{i=1}^3 d_i+1)(k'-\sum_{i=1}^3 d'_i+1)\geq 0.
\]
Since $ L \hookrightarrow \bigoplus_{i=1}^{3} R(-d_i, -d_i')$,
$H_{\textbf{m}}^0(L)=0.$  This implies that 
\begin{equation}
\label{E:vanish}
H_{\textbf{m}}^0(L/K)_{k,k'}=0 \text{ when } (k-\sum_{i=1}^3 d_i+1)(k'-\sum_{i=1}^3 d'_i+1)\geq 0.
\end{equation}
Clearly,  $\mbf{m}^{\alpha} L/K =0$ for some $\alpha$ and by proposition (\ref{P:supp}) and the remark
following, this implies that $\Supp(L/K) \subset \V(\mbf{m})$. It 
is well-known that this last condition implies that   
$H^0_{\mbf{m}}(L/K)=L/K$, and from equation (\ref{E:vanish}) 
this gives $L_{k,k'}=K_{k,k'}$ when $(k-\sum_{i=1}^3 d_i+1)(k'-\sum_{i=1}^3 d'_i+1)\geq 0$.  
This means $r=r_{k,k'} \in K_{k,k'}$
and $K_{k,k'}=K^{\sat}_{k,k'}$ for those same indices, which was to be proved.
\end{proof}

\begin{thm}
\label{T:HZ}
Let $\mc{O}_p$ be the local ring of a point $p$ in $\proj{1}\times
\proj{1}$, and let $\mc{I}_p \subset \mc{O}_p$ be a codimension
two ideal.  Then
\[ dim_k \mc{I}_p/ \mc{I}^2_p \geq 2dim_k \mc{O}_p/ \mc{I}_p.\]
Furthermore, equality holds if and only if $\mc{I}_p$ is a
complete intersection in $\mc{O}_p$.
\end{thm}
\begin{proof}
See \cite{HZ}.
\end{proof}

\begin{thm}
\label{T:LCI}
If $I=\langle f_1, f_2, f_3 \rangle \subset R$ has codimension
two, then $K^{\sat}=V$ if and only if $I$ is a local complete
intersection.
\end{thm}

\begin{proof}
Since $K^{\sat}=V \Leftrightarrow H(K^{\sat})=H(V)$, we will compute both
$H(K^{\sat}), H(V)$.
\par
\[H(K)=H(K^{\sat})= \sum_{i<j}H(R(-d_i-d_j, -d_i'-d_j'))-H(R(-\sum_{i=1}^3 d_i,
-\sum_{i=1}^3 d_i')).\]
\begin{eqnarray*}
H(K^{\sat})(k,k')&=&
\sum_{i<j}(k-d_i-d_j+1)(k'-d'_i-d'_j+1)-(k-\sum_{i=1}^3 d_i +1)
                    (k'-\sum_{i-1}^3 d'_i+1)\\
&=&\sum_{i=1}^3 (k-d_i+1)(k'-d'_i+1)-(k+1)(k'+1)  
\end{eqnarray*}

Now consider $H(V)$. Since $V=S \cap \bigoplus_{i=1}^{3}
I^{\sat}(-d_i, -d_i')$,  and 
\[
0\rightarrow S \rightarrow
\bigoplus_{i=1}^3 R(-d_i, -d'_i) \rightarrow I \rightarrow 0,
\]
we
will have the exact sequence:
\[ 
0 \rightarrow V \rightarrow \bigoplus_{i=1}^{3} I^{\sat}(-d_i,
-d_i') \rightarrow I I^{\sat} \rightarrow 0 . 
\]
 Then
\[ 
H(V)=\sum_{i=1}^3 H(I^{\sat}(-d_i, -d_i'))- H(I I^{\sat}).
\]
Since $ 0 \rightarrow I^{\sat} \rightarrow R \rightarrow
R/I^{\sat} \rightarrow 0$ and $V(I^{\sat})=V(I)=Z$ is
zero-dimensional, we have
\[ 
H(I^{\sat})=H(R)-H(R/I^{\sat})=H(R)-\deg(Z)
\]
Therefore
\[ 
H(V)=\sum_{i=1}^3 H(R(-d_i, -d_i')) -3 \deg(Z) -H(II^{\sat}).
\] 
Note that $I^2, I I^{sat}$ have the same saturation.
To see this, it is enough to show $I I^{\sat} \subset (I^2)^{\sat}$, since $I^2
\subset I I^{\sat}$. Let $f \in I I^{\sat}$, so $f=\sum_{i=1}^k f_i
g_i$ with $f_i \in I, g_i \in I^{\sat}$.  But there exist an $n$
such that $\langle st, sv, ut, uv \rangle ^n g_i \subset I$ for
all $i$, therefore $f \in (I^2)^{\sat}$. Thus from remark
(\ref{R:sat}), we get
$H(I^2)=H(II^{\sat})$. Now
\[ 
H(V)=\sum_{i=1}^3 H(R(-d_i, -d_i')) -3 \deg(Z) -H(I^2).
\]
The exact sequence
\[ 0 \rightarrow I^2 \rightarrow R \rightarrow R/I^2 \rightarrow
0\] and \[ 0 \rightarrow I/I^2 \rightarrow R/I^2 \rightarrow R/I
\rightarrow 0 \] gives
\[ H(I^2)=H(R)-H(R/I^2)=H(R)-H(I/I^2)-\deg(Z).\]
Therefore, we have 
\begin{eqnarray*}
H(V)(k,k')&=&\sum_{i=1}^3H(R(-d_i, -d'_i))(k,k')-H(R)(k,k')-2 \deg (Z)+H(I/I^2)(k,k') \\
 &=& \sum_{i=1}^3 (k-d_i+1)(k'-d'_i+1)-(k+1)(k'+1)- 2\deg(Z) +H(I/I^2)(k,k')
\end{eqnarray*}

Comparing $H(K^{\sat})$ and $H(V)$, we see
\[H(K^{\sat})=H(V) \iff   H(I/I^2)=2 \deg(Z).\]

If $\mc{I}$ is the ideal sheaf of $Z$, then
\[ \deg(Z) =\dim_k H^0(Z, \mc{O}_Z) =\dim_k H^0(\proj{1} \times
\proj{1}, \mc{O}_{\proj{1} \times \proj{1}}/\mc{I}) =\sum_{p\in Z}
\dim_k \mc{O}_p/ \mc{I}_p,\] where $\mc{O}_p, \mc{I}_p$ is the
localization at $p\in Z$.  Since $\mc{I}/\mc{I}^2$ has zero dimensional
support, we have
\[ H(I/I^2)=\dim_k H^0( \proj{1}\times \proj{1}, \mc{I}/\mc{I}^2)=
\sum_{p\in Z} \dim_k \mc{I}_p/\mc{I}^2_p.\] By theorem
(\ref{T:HZ}), we know that \[ \dim_k \mc{I}_p/\mc{I}^2_p \geq 2 \dim_k
\mc{O}_p/ \mc{I}_p \] for every $p \in Z$ with equality holds 
if and only if $\mc{I}_p$ is LCI.  Therefore, we have
\[ H(I/I^2)=2\deg(Z) \Leftrightarrow \dim_k \mc{I}_p/\mc{I}^2_p = 2 \dim_k \mc{O}_p/
\mc{I}_p,  \  \ \forall p \in Z, \]and we conclude that
\[ H(I/I^2)=2\deg(Z) \Leftrightarrow I \text{ is  LCI. } \]
\end{proof}

\begin{cor}\label{C:lcikoz}
If $I=\langle f_1, f_2, f_3 \rangle \subset R$ has codimension
two, and the bidegree of $f_i$ is $(d_i, d'_i)$,  
then $K_{k,k'}=V_{k,k'}\subset \bigoplus_{i=1}^3 R(-d_i, -d'_i)_{k,k'}$ 
when $(k-\sum_{i=1}^3 d_i+1)(k'-\sum_{i=1}^3 d'_i+1)\geq 0$, 
if and only if $I$ is a local complete intersection.
\end{cor}
\begin{proof}
$(\Leftarrow )$ If $I$ is LCI, then $K^{\sat}=V$ by the theorem. 
But $K_{k,k'}=K^{\sat}_{k,k'}$ when $(k-\sum_{i=1}^3 d_i-1)(k'-\sum_{i=1}^3 d'_i+1)\geq 0$, by Lemma \ref{L:sat}.

$(\Rightarrow )$ If $K_{k,k'}=V_{k,k'}$ when 
$(k-\sum_{i=1}^3 d_i-1)(k'-\sum_{i=1}^3 d'_i-1)\geq 0$, 
then $H(K)=H(K^{\sat})=H(V)$.  By Lemma \ref{L:sat}, we have $K^{\sat}=V$.  
By Theorem \ref{T:LCI}, we have $I$ is LCI.
\end{proof}

\begin{rem}
If $d_i=m$ and $d'_i=n$ for $i=1,2,3$, then Corollary \ref{C:lcikoz} implies Proposition \ref{P:vankoz}.  It is because $K_{3m-1,3n-1} \subset \bigoplus_{i=1}^3 R(-m,-n)$ has pure degree $(2m-1,2n-1)$ by Remark \ref{R:puredegree}.
\end{rem}

\begin{exa}
The following example shows that not all syzygies vanishing at the base
point are Koszul syzygies:

$I=\langle s^2v^2, u^2t^2, s^2t^2 \rangle$.  The only base point of ideal
$I$ is $p=(0:1;0:1)$, and $I_p=\langle s^2, t^2 \rangle$.  The base point
is a local complete intersection.  Consider the syzygies $(sut^4v, 0,
-sut^2v^3)$ and $(0, s^4utv, -s^2u^3tv)$ of bidegree $(2,5)$ and $(5,2)$
respectively.  By definition, they vanish on the base point, 
since $\mbf{m}(sut^4v, 0, sut^2v^3) \subset I$, 
$\mbf{m}(0, s^4utv, -s^2u^3tv)\subset I$.
But neither one of them is a Koszul syzygy, 
since the Koszul syzygies are generated by
\[(u^2t^2, -s^2v^2, 0), (s^2t^2, 0, -s^2v^2), (0, s^2t^2, -u^2t^2).\]
\end{exa}

\end{document}